\documentclass{amsart}

\newtheorem{theorem}{Theorem}
\newtheorem{proposition}[theorem]{Proposition}

\newtheorem{corollary}[theorem]{Corollary}

\newtheorem{definition}[theorem]{Definition}

\theoremstyle{remark}
\newtheorem{remark}[theorem]{\bf {Remark}}

\newtheorem{example}[theorem]{\bf {Example}}

\numberwithin{equation}{section}

\usepackage{color}
\usepackage[top=2.5cm, bottom=2cm, left=2.3cm, right=2.5cm]{geometry}

\newcommand{\Nset}{\mathbb{N}}% Si no te funciona, escribe \newcommand{\Nset}{\mathbb{N}}

\newcommand{\Rset}{\mathbb{R}}

\begin{document}

\title{The Samuel realcompactification }
\author {M. Isabel Garrido}
\address{Instituto de Matem\'{a}tica Interdisciplinar (IMI), Departamento de Geometr\'{i}a y Topolog\'{i}a, Universidad Complutense de Madrid, 28040 Madrid, Spain}
\email {maigarri@mat.ucm.es}

\author {Ana S. Mero\~{n}o}
\address {Departamento de An\'{a}lisis Matem\'{a}tico, Universidad Complutense de Madrid, 28040 Madrid, Spain}
\email {anasoledadmerono@ucm.es}

\subjclass [2010] {Primary 54D60, 54E15; Secondary 54D35, 54A20, 54A25}
\keywords{Uniform space, realcompactification,  real-valued uniformly continuous function, Samuel realcompactification, Cauchy filter, Bourbaki-Cauchy filter, Bourbaki-completeness. \\ \indent Partially supported by MINECO Project-MTM2015-65825-P (Spain)}

\begin {abstract} For a uniform space $(X,\mu)$, we introduce a  realcompactification of $X$ by means of the family $U_\mu(X)$ of all the real-valued uniformly continuous functions, in the same way that the known Samuel compactification is given by $U^*_\mu(X)$ the set of all the bounded functions in $U_\mu(X)$.  We will call it ``the Samuel realcompactification'' by several resemblances to the Samuel compactification. In this note, we present different ways to construct such realcompactification as well as we  study the corresponding problem of knowing when a uniform space is Samuel realcompact, that is, it coincides with its Samuel
realcompactification. At this respect we obtain as main result a theorem of Kat\v{e}tov-Shirota type, by means of a new  property of completeness recently introduced by the authors, called Bourbaki-completeness.
\end{abstract}

\maketitle

\begin{section}{Introduction}

A realcompactification of a Tychonoff space $X$ is a realcompact space $Y$ in which $X$ is densely embedded. For instance,  the well-known Hewitt-Nachbin realcompactification   $\upsilon X$. Recall that $\upsilon X$ is  characterized as the smallest realcompactification of $X$ (in the usual order on the family of all the realcompactifications) such that every real-valued continuous function $f\in C(X)$ can be continuously extended to it \cite{gillman}.  In the frame of  uniform spaces, since we can also  consider  $U_{\mu}(X)$, the set of all the real-valued uniformly continuous functions on  the uniform space $(X,\mu)$, it is natural to ask what is the smallest realcompactification of $X$ such that every  function $f\in U_{\mu}(X)$ can be continuously extended to it. Here, following the ideas of \cite{garrido1}, we will  introduce this realcompactification that we denote by  $H(U_{\mu}(X))$ because in fact it can be represented as the set of all the real homomorphisms on the unital vector-lattice $U_{\mu}(X)$. Moreover, we will see that it also coincides with the completion of the uniform space $(X,w_{U_\mu(X)})$, where $w_{U_\mu (X)}$ denotes the weak uniformity on $X$ generated by $U_\mu (X)$ (Theorem \ref{completion}).

\smallskip

We will call $H(U_{\mu}(X))$  the \textit{Samuel realcompactification} of $(X,\mu)$ in likeness to the Samuel compactification $s_{\mu} X$.  Recall that the Samuel compactification of a uniform space $(X,\mu)$, also known as the Smirnov compactification, is the smallest (real)compactification of $X$ such that every bounded real-valued uniformly continuous function $f\in U^{*}_{\mu}(X)$ can be continuously extended to it (\cite{samuel}). The resemblance between the Samuel realcompactification and the Samuel compactification, is not only due to their characterization as the smallest realcompactification, respectively compactification, such that every real-valued, respectively  bounded, uniformly continuous function can be continuously extended to it. In fact, as well as every compactification of a Tychonoff space $X$ can be considered as a Samuel compactification for some compatible (we will say that a uniformity is compatible if it defines the same topology on  $X$)  totally bounded (precompact) uniformity on $X$ \cite{gal} (see also \cite{google}), every realcompactification of $X$ can be considered as a Samuel realcompactification for some compatible uniformity on it, as we will explain forward (Theorem \ref{all.samuel}).

\smallskip
As a Tychonoff space $X$ is realcompact whenever $X=\upsilon X$, then  we can similarly define for a uniform space $(X, \mu)$ to be {\it Samuel realcompact} if $X=H(U_{\mu}(X))$. This paper is mainly devoted to give a uniform analogue to the well-known Theorem of Kat\v{e}tov-Shirota \cite{shirota} (see also \cite{gillman}) about the realcompactness of Tychonoff spaces. Recall that this  classical  theorem states that \textit{a topological Tychonoff space $X$ is realcompact if and only if is completely uniformizable and  every closed discrete subspace has non-measurable cardinal}.  In this line, we will obtain our analogous result  that asserts  that \textit{a uniform space is Samuel realcompact if and only if it is Bourbaki-complete and every uniformly closed subspace has non-measurable cardinal} (Theorem \ref{shirota-like}). The property of Bourbaki-completeness was introduced and well-studied by the authors in \cite{merono1} (see also \cite {merono2} and \cite{merono4}), and it is a uniform property stronger than usual completeness. We say that a uniform space is  \textit{Bourbaki-complete} if every \textit{Bourbaki-Cauchy} filter, defined later, clusters.

It turns out that Bourbaki-Cauchy filters are a useful tool to study Samuel realcompactness. In fact, the Samuel realcompactification, as the completion of $(X,w_{U_{\mu}(X)})$, is always a Bourbaki-complete uniform space (Theorem \ref{H(U) is Bourbaki-complete}). Furthermore, assuming that the uniform partitions of $(X,\mu)$ have non-measurable cardinality, the Cauchy filters $(X,w_{U_{\mu}(X)})$ are precisely Bourbaki-Cauchy filters of $(X,\mu)$   (Proposition \ref{measurable1}), and then it follows that $(X,w_{U_{\mu}(X)})$ is complete if and only of $(X,\mu)$ is Bourbaki-complete (Theorem \ref{measurable2}). And this will be the key not only for proving our main result but in order to find another description of the Samuel realcompactification $H(U_{\mu}(X))$. Namely, $H(U_{\mu}(X))$ is the subset of $s_\mu X$ formed by the cluster points of the Bourbaki-Cauchy filters in $(X,\mu)$ (Theorem \ref{measurable3}).

\smallskip

Very recently in \cite{merono3}, we have proved the same result of Samuel realcompactness but in the frame of metric spaces. The techniques that we  used  are not the same as those of now  since we worked there with a family of metrics defined on the space being  uniformly equivalent to the initial one. Moreover, in this line,  Hu\v{s}ek and Pulgar\'in gave in \cite{husek} a kind of uniform  Kat{e}tov-Shirota result for the case of the so-called uniformly $0$-dimensional uniform spaces. We will see how it can be easily derived now.

\smallskip
In the last section we are going to see that the Samuel realcompactification $H(U_{\mu}(X))$ also  admits a characterization by means of some compatible uniformity whose uniform covers are countable, which reminds in some sense, as we will show, to some of the classical results by Shirota \cite{shirota}.
\end{section}
\medskip

\begin{section}{Preliminaries results  about realcompactifications}

We start this section with some basic facts about realcompactifications, that can be found  mainly in \cite{garrido1}. Thus, the classical  way of generating realcompactifications of a Tychonoff space $X$ is the following. First, we take a family $\mathcal L$ of real-valued continuous functions, that we suppose having the algebraic structure of unital vector lattice and separating points from closed sets of $X$. Then, we embed (homeomorphically) $X$ into the product space of real lines $\Rset ^{\mathcal{L}}$, through the evaluation map $$e: X\rightarrow \Rset ^{\mathcal{L}}$$ $$\hspace{2.7cm} x\rightsquigarrow e(x)=(f(x))_{f\in \mathcal{L}}.$$ Next, we take the closure of $X$ in $\Rset ^{\mathcal{L}}$.  We will denote this closure by $H(\mathcal{L})$ because it is exactly the set of all the real unital homomorphisms on $\mathcal{L}$.  Since $H(\mathcal{L})$ is closed in $\Rset ^{\mathcal{L}}$ then it is in fact a realcompactification of $X$. If we just take the bounded functions $\mathcal{L}^{*}=\mathcal{L}\cap C^{*}(X)$  in $\mathcal{L}$ (where $C^{*}(X)$ is the family of bounded real-valued continuous functions) we get that $H(\mathcal{L}^{*})$ is now a compactification of $X$.

Likewise compactifications, we can consider a partial order $\leq$ on the set $\mathfrak{R}(X)$ of all the realcompactifications of $X$. Namely, for two realcompactifications $\alpha _1 X$ and $\alpha _2 X$, we write $\alpha_1 X\leq \alpha_2 X$ whenever there is a continuous mapping $h:\alpha_2X\to \alpha_1X$ leaving $X$ pointwise fixed. We say that $\alpha_1 X$ and $\alpha_2X$  are equivalent whenever $\alpha_1 X\leq \alpha_2 X$ and $\alpha_2 X\leq \alpha_1 X$, and this implies the existence of a homeomorphism between $\alpha_1 X$ and $\alpha_2 X$ leaving $X$ pointwise fixed.

\smallskip

In particular, $H(\mathcal L)$ (resp. $H(\mathcal L^*)$) is characterized  (up to equivalence) as the smallest realcompactification (resp. compactification) of $X$ such that every function $f\in \mathcal{L}$ can be continuously extended to it. Besides, it is easy to see that  $H(\mathcal L)$ can be considered as a topological subspace of $H(\mathcal L^*)$. Thus, we can write $$(\clubsuit) \hspace{.8cm} X \subset H(\mathcal L)\subset H(\mathcal L^*).$$

If we consider in the family $\mathfrak{R}(X)$ of all the realcompactifications of $X$, the above defined partial order $\leq$, then $(\mathfrak{R}(X), \leq)$ is a complete upper semi-lattice where the largest element is exactly the Hewitt-Nachbin realcompactification $H(C(X))=\upsilon X$. Recall that, in the case of the compactifications $\mathfrak{K}(X)$ of the space $X$, then
$(\mathfrak{K}(X), \leq)$ is also a complete upper semi-lattice where, the largest element is now the Stone-\v{C}ech compactification $H(C^*(X))=\beta X$.

\smallskip

In particular, in the frame of uniform spaces, when we deal with the the Samuel realcompactification and the Samuel compactification defined in the introduction, we have that they are respectively $H(U_{\mu}(X))$ and $s_{\mu} X = H(U^*_{\mu}(X))$, and also that $ H(U_{\mu}(X))\leq \upsilon X$ and $s_{\mu} X \leq \beta X.$
Moreover, by $(\clubsuit)$ we know that $X\subset \upsilon X\subset \beta X  \text{ and } X \subset H(U_{\mu}(X))\subset s_{\mu}X. $
\smallskip

In both cases, for realcompactifications and compactifications, it is known that $( \mathfrak{R}(X), \leq)$, (respectively $( \mathfrak{K}(X), \leq)$) is a complete lattice if and only if $X$ is locally compact. In that case, the smallest element in both lattices is the Alexandroff or the one-point compactification of $X$ which is generated by all the real-valued  functions which are  constant at infinity \cite{mack}.
\smallskip

We will say that a Tychonoff space $X$ is $\mathcal{L}$-realcompact, for the unital vector lattice  of real-valued continuous functions $\mathcal{L}$, if $X=H(\mathcal{L})$. Thus, a space $X$ is realcompact if and only if $X=\upsilon X$ and $X$ is compact if and only if $X=\beta X$. Clearly every $\mathcal{L}$- realcompact space is realcompact. We will see later that, in general, when we are considering different lattices $\mathcal{L}$ and $\mathcal{L}'$, if $X$ is $\mathcal{L}$-realcompact then it is not necessarily true that $X$ is $\mathcal{L}'$-realcompact. For instance, there are realcompact spaces which are not realcompact for other lattices $\mathcal{L}$ different from $C(X)$. However, when a space is compact then it is $\mathcal{L}^{*}$-compact for any lattice $\mathcal{L}$.
\smallskip

On the other hand, for a unital vector lattice  $\mathcal{L}\subset C(X)$ we can consider $w _{\mathcal{L}}$ the weak uniformity in $X$ which is the weakest uniformity making each function in $\mathcal{L}$ uniformly continuous \cite{willard}. When $\mathcal{L}$ separates points and closed sets in $X$, then this Hausdorff uniformity is compatible with the topology of $X$. If we endowed $X$ with the weak uniformity $w_{\mathcal{L}}$ and $\Rset ^{\mathcal{L}}$ with the product uniformity then, the evaluation map $e:X\rightarrow \Rset ^{\mathcal{L}}$ is uniformly  continuous and the inverse map $$e^{-1}: e(X)\rightarrow X$$ is also uniformly continuous. Thus, $X$ is uniformly embedded in $\Rset ^{\mathcal{L}}$. Since $\Rset ^{\mathcal{L}}$ endowed with the product uniformity is a complete uniform space then $H(\mathcal{L})$, being the closure of $X$ in $\Rset ^{\mathcal{L}}$, must be the completion of $(X, w_{\mathcal{L}})$ by uniqueness of the completion.  We can summarize all of this  as follows.

\begin{theorem} \label{completion} The realcompactification $H(\mathcal{L})$, where $\mathcal{L}\subset C(X)$ is a unital vector lattice  separating points and closed sets in $X$, is homeomorphic to the completion of the uniform space $(X, w_{\mathcal{L}})$ where $w_{\mathcal{L}}$ is the weak uniformity generated  by $\mathcal{L}$.
\end{theorem}

If we apply this result to the above defined realcompactifications we get that $\upsilon X$ is homeomorphic to the completion of $(X, w_{C(X)})$, $\beta X$ to the completion of $(X, w_{C^{*}(X)})$, $H(U_{\mu}(X))$  to the completion of $(X, w_{U_{\mu}(X)})$  and $s_{\mu}X$ to the completion of $(X, w_{U^{*}_{\mu}(X)})$.
\smallskip

In general, when we have a realcompactification $Y$ of $X$ which is not generated by an explicit lattice $\mathcal{L}$,  if $C(Y)$ is the family of all the real-valued continuous functions on $Y$, then $(Y, w_{ C(Y)})$ is complete \cite{gillman}. Precisely it is the completion of $(X, w_{C(Y)_{|X}})$ where $C(Y)_{|X}\subset C(X)$ are the restrictions to $X$ of the functions in $C(Y)$ \cite{hager-johnson}. Moreover we can describe  $C(Y)_{|X}$ as the continuous real-valued functions on $X$ which preserves Cauchy filters of the weak uniformity $w_{C(Y)_{|X}}$ \cite{borsik}. Therefore since $U_{w_{C(Y)_{|X}}}(X)=C(Y)_{|X}$, because every uniformly continuous functions preserves Cauchy filters (or nets) \cite{borsik}, we have the following result.

\begin{theorem}\label{all.samuel} Let $Y$ be a realcompactification of the Tychonoff space $X$. Then $Y$ is homeomorphic to the Samuel realcompactification of the uniform space $(X, w_{C(Y)_{|X}})$.
\end{theorem}

For instance, the Hewitt-Nachbin realcompactification $\upsilon X$ can be considered the Samuel realcompactification of the uniform space $(X, w_{C(X)})$ \cite{gillman}. Equivalently, $\upsilon X$ is the Samuel realcompactification of the uniform space $(X,u)$ where $u$ is the universal uniformity on $X$, because the family of real-valued uniformly continuous functions on $(X,u)$ is exactly $C(X)$. In the same way, the Stone-\v{C}ech compactification $\beta X$ is the Samuel compactification of the uniform space $(X, w_{C^{*}(X)})$.

\smallskip

Now, we are going to finish this section with some easy results in this topic, that we will use later.

\begin{theorem} \label{eq.Samuel} Let $(X,\mu)$ be a uniform space. For every compatible uniformity $\nu$ on $X$ satisfying that $w_{U_{\mu}(X)} \preceq \nu \preceq \mu$ it is satisfied that $H(U_{\mu}(X))=H(U_{\nu}(X))$. Besides, $s_{\mu} X= s_{\nu} X$.
\begin{proof} It follows easy from the fact that $U_{\mu}(X) =U_{w_{U_{\mu}(X)}} (X)\subset U_{\nu}(X)\subset U_{\mu}(X)$. The same is true for bounded functions.
\end{proof}
\end{theorem}

\begin{theorem} \label{eq.Samuel.complete}Let $(\widetilde{X},\widetilde{\mu})$ be the completion of a uniform space $(X,\mu)$. Then $H(U_{\mu}(X))=H(U_{\widetilde{\mu}}(\widetilde{X}))$ and $s_{\mu}X =s_{\widetilde{\mu}}\widetilde{X}$.
\begin{proof} It is known that the functions in $U_{\mu}(X)$ are exactly the restrictions of the real-valued uniformly continuous functions of the completion $(\widetilde{X}, \widetilde{\mu})$ (see \cite{willard}). Thus, by density, the result follows. As in the previous theorem, the same is true for bounded functions.
\end{proof}
\end{theorem}

\end{section}
\medskip

\begin{section}{Main results}

We star this section recalling the notion of Bourbaki-completeness for uniform spaces that was introduced and well studied by the authors in \cite{merono1}, in the frame of metric spaces, and in \cite{merono2} for some special uniformities.

So let $\mathcal{U}\in \mu$, a uniform cover in $(X,\mu)$. For $U\in \mathcal{U}$ let us write $$St^1(U, \mathcal{U})= St(U, \mathcal{U}):=\bigcup \{V\in\mathcal{U}: V\cap U\neq\emptyset \}$$ and put $$St^{m}(U,\mathcal{U})=St (St ^{m-1}(U, \mathcal{U}), \mathcal{U} ),\text{ }m\geq 2$$ $$St^{\infty}(U,\mathcal{U})=\bigcup _{m=1}^{\infty} St ^{m}(U, \mathcal{U}).$$

\begin{definition} A filter $\mathcal{F}$ in the uniform space $(X,\mu)$ is said to be Bourbaki-Cauchy if for every $\mathcal{U}\in \mu$ there exist  $m\in \Nset$ and $U \in \mathcal{U}$ such that $F\subset St^{m}(U,\mathcal{U})$, for some $F\in \mathcal{F}$ (i.e. $St^{m}(U,\mathcal{U})\in \mathcal{F}$ ). Moreover,  we said that $(X,\mu)$ is  Bourbaki-complete whenever every Bourbaki-Cauchy filter in $X$ clusters.
\end{definition}

It is easy to see that every Cauchy filter is Bourbaki-Cauchy, and therefore every Bourbaki-complete uniform space is complete. In general the reverse implication is not true. For example, every infinite-dimensional Banach space is not Bourbaki-complete with the uniformity given by its norm (see \cite{merono1}). In fact, we know that a normed space is Bourbaki-complete if and only if it has finite dimension \cite{merono1}. More examples can be found in \cite{merono4}.

On the other hand, it is easy to check that a subspace of a Bourbaki-complete uniform space is Bourbaki-complete if and only if it is closed.  Furthermore this uniform property is also productive as next result proves.

\begin{proposition} \label{product} Any nonempty product of uniform spaces is Bourbaki-complete if and only if each factor is Bourbaki-complete.

\begin{proof} Suppose  $\Pi X_i$ is Bourbaki-complete. Since each factor  $X_i$ is (uniformly) homeomorphic to a closed subspace of this product, then it must be Bourbaki-complete, as we have said above.

On the other hand, suppose $X_i$ is Bourbaki-complete for every $i\in I$, and let $\mathcal{F}$ a Bourbaki-Cauchy filter in the product. Take  $\mathcal{H}$  an ultrafilter containing $\mathcal{F}$. Clearly, $\mathcal{H}$ is also Bourbaki-Cauchy and then its projection into $X_i$ will be a Bourbaki-Cauchy ultrafilter, for every $i\in I$. Now, from the Bourbaki-completeness of every factor, this projection must converges to a point in $X_i$. Therefore,  $\mathcal{H}$ also converges to a point in the product, and this means, in particular, that the initial filter $\mathcal{F}$  clusters, as we wanted.
\end{proof}
\end{proposition}

\begin{theorem} \label{H(U) is Bourbaki-complete} Let $(X,\mu)$ be a uniform space. Then its Samuel realcompactification $H(U_{\mu}(X))$, as the completion of $(X, w_{U_\mu(X)})$,  is Bourbaki-complete.
\begin{proof} First note that the completion of $(X, w_{U_\mu(X)})$ is precisely  the closed subspace $H(U_{\mu}(X))\subset\mathbb{R}^{U_\mu(X)}$  with the uniformity inherit by the product space $\mathbb{R}^{U_\mu(X)}$. And  then the proof follows at once from the above Proposition \ref{product}, since $\mathbb{R}$ is Bourbaki-complete.
\end{proof}
\end{theorem}

Next result, that is an easy corollary of the above, is in the line
of those contained in \cite{merono2} comparing both completeness properties for some  special uniformities.

\begin{theorem} \label{w-complete = w-Bourbaki-complete} Let $(X,\mu)$ be a uniform space. Then $(X,w_{U_\mu (X)})$ is complete if and only if it  is Bourbaki-complete.
\end{theorem}

Now we are going to see an interesting relationship between the Bourbaki-Cauchy filters in the space  $(X,\mu)$ and its Samuel realcompactification.

\begin{proposition} \label{clusters} Let $(X,\mu)$ be a uniform space and $Z\subset s_{\mu} X$ the subspace of all the cluster points  of the Bourbaki-Cauchy filters of $(X,\mu)$. Then $X\subset Z\subset H(U_{\mu}(X))$.
\begin{proof} Clearly $X\subset Z$. So, let $\mathcal{F}$ be a Bourbaki-Cauchy filter in $(X, \mu)$ and $\xi\in s_\mu X$ a cluster point of $\mathcal{F}$. Consider the filter $\mathcal{H}$ in $X$ generated by the family $\{F\cap V: F\in \mathcal{F},\, V \text{\, is a neighborhood of \,} \xi \text{\, in\, } s_\mu X  \} $. Then  $\mathcal{H}$ is also Bourbaki-Cauchy in $(X,\mu)$, since it contains $\mathcal{F}$. Taking into account  that the identity mapping $id: (X, \mu) \rightarrow (X,w_{U_{\mu}(X)})$ is uniformly continuous, $\mathcal{H}$ is Bourbaki-Cauchy in $(X,w_{U_{\mu}(X)})$. Now, from Theorem \ref{H(U) is Bourbaki-complete},  this filter must cluster in  $H(U_{\mu}(X))$. And we finish, noting that $\xi$ is the only cluster point of $\mathcal{H}$.
\end{proof}
\end{proposition}

Now, we are going to study the problem of the Samuel realcompactness of a uniform space, and we will see that a Kat\v{e}tov-Shirota type theorem  can be obtained where Bourbaki-completeness will play the role of completeness in the classical one. Recall that every uniform countable cover $\{U_n:n\in \mathbb Z\}\in \mu$ such that $U_n\cap U_m=\emptyset$ whenever $|n-m|>1$ belongs also to the weak uniformity $w_{U_\mu (X)}$ (see \cite{isbell}). This kind of covers were called  {\it linear} covers by Isbell in \cite{isbell}, and  2-{\it finite} covers by Garrido and Montalvo in \cite{garrido-montalvo}.

\begin{proposition} \label{measurable1} Let $(X,\mu)$ be a uniform space. If  every uniform partition of $(X,\mu)$ has non-measurable cardinal then every Cauchy filter of $(X,w_{U_\mu (X)})$ is a Bourbaki-Cauchy filter of $(X,\mu)$.
\begin{proof} Let  $\mathcal{F}$ be a Cauchy filter in $(X,w_{U_\mu (X)})$ and let $\mathcal{U}\in \mu$. Note that the family $\{St^{\infty}(U,\mathcal{U}): U\in \mathcal{U}\}$, can be seen as a (uniform) partition of $X$. Indeed, it is clearly a cover of $X$ and two members of this family are equal or they are disjoint.  Now we are going to prove that, since every uniform partition of $(X,\mu)$ has non-measurable cardinal then there exists a unique $U\in \mathcal{U}$ such that $St^{\infty}(U,\mathcal{U})\in \mathcal{F}$. Observe that we can fix representative elements $U_{i}\in \mathcal{U}$, $i\in I$, such that $X=\bigcup \{St^{\infty}(U_{i},\mathcal{U}):i\in I\}$ and $St^{\infty}(U_{i},\mathcal{U}) \cap St^{\infty}(U_{j},\mathcal{U})=\emptyset$ whenever $i\neq j,$ $i,j\in I$. Next, define $$\mathcal{J}:=\left\{J\subset I :\bigcup \{St^{\infty}(U_{i},\mathcal{U}):i\in J  \} \in  \mathcal{F} \right\}.$$ We want to prove that $\mathcal{J}$ is an ultrafilter in  $I$, with the countable intersection property, where $I$  is endowed with the discrete metric $0-1$. We are going to check only that $\mathcal{J}$ satisfies the maximal property for being an ultrafilter and the  countable intersection property, leaving the remaining properties to the reader.

Let $J\subset I$ and suppose that $J\notin \mathcal{J}$. In order to see that $I-J \in \mathcal{J}$, consider  $$\mathcal{W}=\left\{\bigcup \{St^{\infty}(U_{i},\mathcal{U}):i\in J  \},\, \bigcup \{St^{\infty}(U_{i},\mathcal{U}):i\in I-J  \}\right\}.$$ Since  $\mathcal{W}$ is a uniform finite partition of $X$ then  $\mathcal{W}\in w_{U_\mu (X)}$. Now as $\mathcal{F}$ is a Cauchy filter, then we have that $$\bigcup \{St^{\infty}(U_{i},\mathcal{U}):i\in I-J  \} \in \mathcal{F}.$$

Now we prove that $\mathcal{J}$ satisfies the countable intersection property.  Take $\{J_n: n\in \Nset\}\subset \mathcal{J}$ and suppose without loss of generality it is a strictly decreasing family. If $\bigcap \{J_n: n\in \Nset\}=\emptyset$, then the sets $\{I- J_n: n\in \Nset\}$, form a cover on $I$. Define the family of sets $\mathcal{W}=\{W_n:n\in \Nset\}$ by $$W_n=\bigcup \{ St^{\infty}(U_{i},\mathcal{U}):i \in (I- J_{n}) \cap J_{n-1}\},$$
 where $J_0=I$.  Then, $\mathcal{W}$ is a uniform countable partition of $X$ and  therefore, $\mathcal{W}\in w_{U_{\mu}(X)}$. Since  $\mathcal{F}$ is a Cauchy filter of $  (X,w_{U_\mu (X)})$ then $W_n\in \mathcal{F}$ for some $n\in \Nset$. Since $W_n\subset \bigcup \{St^{\infty}(U_{i},\mathcal{U}):i \in I- J_n\}$ then  $\bigcup \{ St^{\infty}(U_{i},\mathcal{U}):i \in I- J_n\} \in \mathcal{F}$, and therefore  $I-J_n \in \mathcal{J}$. But this is a contradiction because $\mathcal{J}$ is an ultrafilter and $J_n\in \mathcal{J}$. Hence, $\bigcap \{J_n: n\in \Nset\} \neq \emptyset$ and $\mathcal{J}$ has the countable intersection property.
\smallskip

By hypothesis, the uniform partition $\{St^{\infty}(U_{i},\mathcal{U}):i\in I\}$ of $X$,  has non-measurable cardinal, and therefore  the discrete space $I$ is realcompact. Hence the ultrafilter $\mathcal{J}$ must be fixed, that is, there exists a unique $i_0\in I$ such that $i_0\in \bigcap \mathcal{J}$. But this is equivalent to say that there exists a unique $i_0\in I$ such that $St^{\infty}(U_{i_0},\mathcal{U}) \in \mathcal{F}$.

Now we are going to see that $\mathcal{F}$ is Bourbaki-Cauchy in $(X,\mu)$. For that we need to find some $m\in \Nset$ such that $St^{m}(U_{i_0},\mathcal{U}) \in \mathcal{F}.$ So, let us define the cover $\mathcal{A}=\{A_n:n\in \{0\}\cup\Nset\}$ of $X$ as follows:

\begin{center}

$A_0=\bigcup \{St^{\infty}(U_{i},\mathcal{U}):i \in I -\{i_0\}\},$
\smallskip

$A_1=St(U_{i_0},\mathcal{U}),$
\smallskip

$A_2=St^{2}(U_{i_0},\mathcal{U}),$
\smallskip

$A_n=\bigcup\{V\in St^{\infty}(U_{i_0},\mathcal{U}):V\cap St^{n}(U_{i_0},\mathcal{U})\neq \emptyset , V\cap St^{n-2}(U_{i_0},\mathcal{U})= \emptyset\}$ if $n\geq 3$.

\end{center} Then $\mathcal{A}$ is a uniform countable cover satisfying that $A_i \cap A_j =\emptyset $ whenever $|i-j|>1$ and hence, $\mathcal{A}\in w_{U_{\mu}(X)}$. Again, since $\mathcal{F}$ is Cauchy in $(X,w_{U_\mu (X)})$ there exists some $n\in \{0\} \cup \Nset$ such that $A_n \in \mathcal{F}$. But as $A_0 \notin \mathcal{F}$, then $A_n \in \mathcal{F}$ for some $n\in \Nset$. Since $A_n\subset St^{n+1}(U_{i_0},\mathcal{U})$, then $St^{n+1}(U_{i_0},\mathcal{U}) \in \mathcal{F}$ and we have finished.
\end{proof}
\end{proposition}

\begin{theorem} \label{measurable2} Let $(X,\mu)$ be a uniform such that every uniform partition of  $X$ has non-measurable cardinal. Then, $(X,\mu)$ is Bourbaki-complete if and only if $(X,w_{U_\mu (X)})$ is complete.
\begin{proof} First, suppose  that $(X,w_{U_\mu (X)})$ is complete then, from Proposition \ref{w-complete = w-Bourbaki-complete}, it is Bourbaki-complete. Now, as the uniformity $\mu$ is finer than $w_{U_\mu (X)}$, then $(X,\mu)$ will be also Bourbaki-complete.

Reciprocally, suppose $(X,\mu)$ is Bourbaki-complete and let $\mathcal{F}$ be a Cauchy filter in $(X,w_{U_\mu (X)})$. From Proposition \ref{measurable1}, we have that $\mathcal{F}$ is  a Bourbaki-Cauchy filter in $(X,\mu)$, and hence it clusters. Since any Cauchy filter with a cluster point converges, it follows that $(X,w_{U_\mu (X)})$ is complete.
\end{proof}
\end{theorem}

We have already all the ingredients in order to establish our main result in this section.

\begin{theorem}\label{shirota-like} {\rm ({\sc Kat\v{e}tov-Shirota type theorem)}} Let $(X,\mu)$ be a uniform space. Then $(X,\mu)$ is Samuel realcompact if and only $(X,\mu)$ is Bourbaki-complete and every  uniform discrete subspace  of  $X$ has non-measurable cardinal.
\begin{proof} If $(X,\mu)$ is Samuel realcompact then it is realcompact and hence every discrete closed subspace has non-measurable cardinal \cite{gillman}. In particular,  every uniformly discrete subspace must have non-measurable cardinal since it is in addition closed. On the other hand, from the identity $X=H(U_\mu(X))$, it follows that  $(X, w_{U_\mu (X)})$ coincides with its completion and then $(X, w_{U_\mu (X)})$ is complete, or equivalently,  from Theorem \ref{measurable2}, $(X,\mu)$  is Bourbaki-complete. Note that we can apply that theorem because every uniform partition produces a  uniformly discrete subspace with the same cardinal.

Conversely, again from Theorem \ref{measurable2}, the Bourbaki-completeness of $(X,\mu)$ together with the property of non-measurability of the corresponding cardinals, imply  the completeness of $(X,w_{U_\mu (X)})$. Therefore  $X=H(U_\mu(X))$, and this means that  $(X,\mu)$ is Samuel realcompact.
\end{proof}
\end{theorem}

\begin{remark} Note that the condition of non-measurable cardinality can not be deleted in the three previous results. Indeed, if $X$ is a set with a measurable cardinal endowed with the uniformity $\mu$ given by the 0-1 metric, then the corresponding uniform (metric) discrete space $(X,\mu)$ is Bourbaki-complete (\cite{merono1}), but not Samuel realcompact since it is not even realcompact.  In this case $X\neq\upsilon X =  H(U\mu(X))$, $(X,w_{U_\mu (X)})$ is not complete, and there exist Cauchy filters in $(X,w_{U_\mu (X)})$ that are not Bourbaki-Cauchy filters in $(X,\mu)$.
\end{remark}

On the other hand, it is clear that in Theorem \ref{shirota-like} only the non-measurability of the uniform partitions are needed and not the (stronger) condition of non-measurability of the uniform discrete subspaces. Next result characterizes precisely the property for a uniform space to have every uniform partition with a non-measurable cardinal. But in order to establish this we need to recall some property of the Samuel compactification (see \cite{engelkingbook}). Namely, for subsets $A$ and $B$ of a uniform space $(X,\mu)$, we have that  $$cl_{s_{\mu}X} A \cap cl_{s_{\mu}X} B \neq \emptyset \Longleftrightarrow  St(A,\mathcal{U})\cap St(B,\mathcal{U})\neq \emptyset,  \text{ for all\,\,} \mathcal{U}\in \mu.$$

\begin{theorem} \label{measurable3}  Let $(X,\mu)$ be a uniform space and $Z$ the subspace of  $s_{\mu}X$ of  all the cluster points  of the Bourbaki-Cauchy filters of $(X,\mu)$. The following statements are equivalent:
\begin{enumerate}

\item every uniform partition of $(X,\mu)$ has non-measurable cardinal;

\item $Z=H(U_{\mu}(X))$;

\item $Z$ is realcompact.

\end{enumerate}
\begin{proof} $(1)\Rightarrow (2)$ By Proposition \ref{clusters}, we have  $Z\subset H(U_{\mu}(X))$. The other inclusion follows from Proposition \ref{measurable1}.

$(2)\Rightarrow (3)$ This implication is trivial.

$(3)\Rightarrow (1)$. We are going to see how every uniform partition of $(X,\mu)$ determines a closed discrete subspace of $Z$. Indeed, let $\mathcal{V}$ a uniform partition of $X$. Then for some (open) uniform cover $\mathcal{U}\in \mu$ , $St (V,\mathcal{U})\cap St (V',\mathcal{U})=\emptyset$, for every $V \neq V'$, $V, V' \in \mathcal{V}$. In particular, for $V\in\mathcal{V}$, we have that $V=St(V,\mathcal{U})$ and this implies that $V$ is open and $V= St^m(V, \mathcal{U})$, for every $m\in \mathbb N$. Moreover,  if $\mathcal{F}$ is a Bourbaki-Cauchy filter of $(X,\mu)$ then there is some (unique) $V\in \mathcal{V}$ such that $V\in \mathcal{F}$. Thus if $\xi\in s_{\mu}X$ is a cluster point of $\mathcal{F}$ then $\xi \in cl_{s_{\mu}X} V$. Hence, $Z\subset \bigcup \{cl _{s_{\mu} X} V: V\in \mathcal{V}\}$, where  $cl_{s_\mu X} {V}\cap cl_{s_{\mu }X}V'=\emptyset$,  for every $V, V' \in \mathcal{V}$ with $V\neq V'$ (see the equivalence before the theorem).

Now, for every $V\in \mathcal{V}$ take a representative point $x_V\in V$. Then $D=\{x_V:V\in\mathcal{V}\}$ is a discrete subspace of $s_\mu X$. Indeed, every $cl_{s_\mu X} {V}$ is a neighborhood of every point belonging to $V$, because if $x\in V=\widetilde V\cap X$ for some open set $\widetilde V$ of $s_\mu X$, and then $x\in \widetilde V\subset cl_{s_\mu X} {\widetilde V}=cl_{s_\mu X} {V}$. Finally in order to see that $D$ is closed in $Z$, suppose $\xi\in Z\setminus D$, then there exists a unique $V\in \mathcal{V}$
such that $\xi\in cl _{s_{\mu} X} V$. Again as $St(V,\mathcal{U})\cap St(X\setminus V,\mathcal{U}) = \emptyset$, we have that $cl_{s_\mu X} {V}\cap cl_{s_{\mu }X}(X\setminus V)=\emptyset$. This means that $\xi \notin cl_{s_{\mu} X }D$.

We finish, since if $Z$ is realcompact then $D$ has non-measurable cardinal. Therefore the uniform partition $\mathcal{V}$ has non-measurable many elements.
\end{proof}
\end{theorem}

Note that last result provides another interesting  description of the Samuel realcompactification for those uniform spaces fulfilling condition (1). Namely, the Samuel realcompactification of these spaces are formed by the cluster points in the Samuel compactification of their Bourbaki-Cauchy filters. Note that there are many spaces with this property,  for instance, connected spaces, or more generally uniformly connected, separable,  Lindel\"{o}f, and many other spaces.

Observe also  that this last theorem reminds of the result in \cite[Theorem 15.21]{gillman} that asserts that the completion of a uniform space is realcompact if and only if every uniform closed discrete subspace has non-measurable cardinal. In fact, think of the completion of a uniform space as the set of all the convergence points, equivalently cluster points, in $s_{\mu}X$ of the Cauchy filters in $(X,\mu)$, and take into account that every uniform partition is a uniformly closed subspace.  However, in our case, we cannot do without Bourbaki-completeness and we cannot change uniform partitions by uniformly closed subspaces.

\smallskip

We finish this section recalling a result of Kat\v{e}tov-Shirota type given by Hu\v{s}ek and Pulgar\'{\i}n in the setting of uniform spaces with weak uniformities (\cite{husek}). They define a uniform space $(X, \mu)$ to be {\it  uniformly realcomplete} when $(X, w_{U_\mu X})$ is complete. Thus, in that paper,  it is proved the following result for the particular case of uniformly 0-dimensional spaces, where a uniform space is {\it uniformly 0-dimensional} if it has a base for the uniformity composed by partitions (for instance, every uniformly discrete space).

\begin{theorem} {\rm ({\sc Hu\v{s}ek-Pulgar\'{\i}n} \cite{husek})}  A 0-dimensional uniform space is uniformly realcomplete if and only if it is complete  and it does not have uniformly discrete subsets of measurable cardinality.
\begin{proof} The proof follows easily from Theorem \ref{shirota-like}. Indeed, we only need to take into account that every complete uniformly 0-dimensional space is clearly Bourbaki-complete, and also that for the uniform space $X$ the condition to be  uniformly realcomplete is equivalent to say that  $X$ is Samuel realcompact.
\end{proof}
\end{theorem}

\end{section}
\medskip
\begin{section}{Uniformities with uniform countable covers}

In his famous paper Shirota \cite{shirota} proved that the Hewitt-Nachbin realcompactification $\upsilon X$ of a space $X$ is homoemorphic to the completion of $(X, eu)$ where $eu$ is the uniformity having as a base all the countable normal covers of $X$. Thus he was proving that the completion of $(X, w_{C(X)})$ and the completion of $(X, eu)$ are (topologically) equivalent even if the uniformities $w_{C(X)}$ and $eu$ are not.

\smallskip

Later on, Isbell proved (see  \cite{isbellbook}) a similar result in the frame of the so called \textit{locally fine} uniform spaces (\cite{ginsburg}). Recall that to be locally fine is equivalent to the more intuitive notion of being \textit{subfine}, i.e., to be a uniform subspace of a fine space \cite{pelant}. The result of Isbell states that for a locally fine uniform space $(X,\mu)$ every Cauchy filter in $(X, w_{U_{\mu}(X)})$ is Cauchy in $(X, e\mu)$, where $e\mu$ is the compatible uniformity on $X$ having as a base all the countable covers of $\mu$ (\cite{ginsburg}). Then we can say now that, for locally fine uniform space $(X,\mu)$, the completion of $(X, e\mu)$ is homeomorphic to its Samuel realcompactification $H(U_{\mu}(X))$.

\smallskip

In particular, since every fine space is locally fine, we can apply last result to the fine uniformity $u$ on a Tychonoff space $X$. Thus, we can deduce easily the above result of Shirota because for fine uniform spaces every continuous function is uniformly continuous. Observe that, in general, locally fine uniform spaces do not satisfy that every continuous function is uniformly continuous.
In fact, locally fine uniform spaces are characterized as those uniform spaces satisfying that every uniformly locally uniformly continuous function into a metric space is uniformly continuous \cite{pelant}. Therefore, for locally fine uniform spaces the Hewitt-Nachbin  realcompactification and the Samuel realcompactification are not necessarily equivalent.

\smallskip

In the next example we show that, in general, the completion of $(X, e\mu)$ is not homeomorphic to the Samuel realcompactification $H(U_{\mu}(X))$. To that purpose we need to recall some results and concepts from \cite{merono1}.

Recall that a uniform space $(X,\mu)$ is \textit{Bourbaki-bounded} if for every $\mathcal{U}\in \mu$ there exists $m\in \Nset$ and finitely many $U_i \in \mathcal{U}$, $i=1,...,k$ such that $X=\bigcup _{i=1}^{k} St^m(U_i,\mathcal{U})$. Observe that it was precisely the notion of Bourbaki-bounded that originated the study of Bourbaki-Cauchy filters, nets and sequences, as well as  Bourbaki-completeness (see \cite{merono1} and \cite{merono2}). In particular in \cite{merono2}  it is proved that \textit{a uniform space is compact if and only if it is Bourbaki-bounded and Bourbaki-complete}.

\begin{example} Let $(X,d)$ be any complete metric space which is also Lindel\"{o}f and Bourbaki-bounded for the metric uniformity, but not compact.  For instance the metric hedgehog $H(\aleph _0)$ of countable many spines, or any closed and bounded subset of the classical Hilbert space $\ell ^{2}$. Then the metric uniformity given by $d$ coincides with the uniformity $e\mu _d$ by the Lindel\"of property. Therefore $(X, e\mu _d)$ is complete. However, by Theorem \ref{measurable2},   $(X,w_{U_d(X)}(X))$ is not complete because $(X,d)$ fails to be Bourbaki-complete. Otherwise, since $(X,d)$ is Bourbaki-bounded $X$ would be compact which is false.
\end{example}

In the previous example, the completion of $(X,e\mu)$ is not homeomorphic to the Samuel realcompactification $H(U_{\mu}(X))$ because completeness is not enough. In fact we need Bourbaki-completeness as it is shown by the following two results.

\begin{theorem} \label{countable1} Let $(X,\mu)$ be a uniform space and $Y$ the subspace of $s_{\mu} X$ of all the cluster points of the Bourbaki-Cauchy filters of $(X,e\mu)$. Then $H(U_{\mu}(X))=Y$.
\begin{proof} By Theorem \ref{eq.Samuel}, it is clear that $H(U_{\mu}(X))=H(U_{e\mu}(X))$. Since it is not difficult to see that every uniform partition of $(X,e\mu)$ is countable the result follows from  Theorem \ref{measurable3}.
\end{proof}
\end{theorem}

\begin{corollary}\label{countable2} A uniform space $(X,\mu)$ is Samuel realcompact if and only if $(X,e \mu)$ is Bourbaki-complete.
\end{corollary}

\end{section}

\end{document}